%% file: go4rdec10.tex
\documentclass[12pt]{article}
\usepackage[english]{babel}
\usepackage[english]{babel}
\usepackage{lineno} 
\usepackage{graphicx,multicol}
\usepackage{epic,eepic,epsfig}
\usepackage{amssymb}
\usepackage{amsmath,amsthm}
\usepackage{color}
\usepackage{tikz}

\usepackage{float}

\usepackage{eso-pic} 

\newcommand{\pf}{{\bf Proof: }}

\newcommand{\2}{\vspace{2mm}}

\theoremstyle{plain}
\newtheorem{theorem}{Theorem}
\newtheorem{proposition}[theorem]{Proposition}

\newtheorem{corollary}[theorem]{Corollary}
\newtheorem{lemma}[theorem]{Lemma}

\theoremstyle{definition}

\newtheorem{problem}[theorem]{Problem}

\newtheorem{definition}[theorem]{Definition}

\begin{document}
\bibliographystyle{plain}
\title{Good acyclic orientations of 4-regular 4-connected graphs\thanks{Research supported by the Independent Research Fond Denmark under grant number DFF 7014-00037B.}}
\author{J. Bang-Jensen\thanks{Department of Mathematics and Computer Science, University of Southern Denmark, Odense, Denmark (email: jbj@imada.sdu.dk)}
    \and M. Kriesell\thanks{Department of Mathematics, Technische Universit\"at Ilmenau, Germany (email: matthias.kriesell@tu-ilmenau.de)}}
  \maketitle

  \begin{abstract}
    We study graphs which admit an acyclic orientation that 
contains an  out-branching and in-branching which are arc-disjoint (such an orientation 
is called {\bf good}).
A {\bf 2T-graph} is a graph whose edge set can be decomposed into two edge-disjoint
spanning trees. Clearly a graph has a good orientation if and only if it contains a spanning 2T-graph with a good orientation, implying that 2T-graphs play a central role. It is a well-known result due to Tutte and Nash-Williams, respectively, that every 4-edge-connected graph contains a spanning 2T-graph.
Vertex-minimal 2T-graphs with at least two vertices, also 
known as {\bf generic circuits}, play an important role in rigidity theory for 
graphs. It was shown in \cite{bangGOpaper}  that every generic circuit has a good orientation.
Using this,  several results on good orientations of 2T-graphs were obtained in \cite{bangGOpaper}.
It  is an open problem whether there exist a polynomial algorithm for deciding whether a given 2T-graph has a good orientation. In \cite{bangGOpaper} complex constructions of 2T-graphs with no good orientation were given, indicating that the problem might be very difficult.
In this paper we focus on so-called {\bf quartics} which are 2T-graphs where every vertex has degree 3 or 4. We identify  a sufficient condition for a quartic to have a good orientation, give a polynomial algorithm to recognize quartics satisfying the condition and a polynomial algorithm to produce such an orientation when this
condition is met. As a consequence of these results we prove that every 4-regular and 4-connected graph has a good orientation.
We also show that every graph on $n\geq 8$ vertices and of minimum degree at least $\lfloor{}n/2\rfloor$ has a good orientation. Finally we pose a number of open problems.\\

\noindent{}{\bf Keywords:} acyclic orientations, edge-disjoint spanning trees, generic circuit, 4-regular graph, 2T-graph, polynomial algorithm, out-branching, in-branching.

  \end{abstract}

  \section{Introduction}

A graph $G=(V,E)$ on at least two vertices is a {\bf 2T-graph} if it has two spanning trees $S,T$ such that $E(S) \cap E(T)=\emptyset$ and $E(S) \cup E(T)=E$,
and it is called a {\bf generic circuit} if, moreover, it does not contain a $2T$-graph as a proper subgraph.
Observe that any 2T-subgraph of a 2T-graph is an induced subgraph. In particular, the generic circuits of a 2T-graph are
induced subgraphs. A $2$-cycle (formed by a pair of parallel edges)
is a generic circuit, which implies that a generic circuit on more than two vertices is simple. Since
deleting a vertex of degree $2$ from a 2T-graph produces a 2T-graph, these larger generic circuits will have minimum degree $3$.

It is well-known that one can decide the existence of a spanning 2T-graph in a graph $G$ in polynomial time
(see e.g. \cite{recski1989}).

Let $G=(V,E)$ be a graph. For a given partition ${\mathfrak F}$ of the vertices of $G$ we denote by $E_G({\mathfrak F})$ the set of edges whose end vertices lie in different sets of $\mathfrak F$. 
The following theorem, due to Nash-Williams and Tutte, characterizes graphs with $k$ edge-disjoint spanning trees.

\begin{theorem} \cite{nashwilliamsJLMS39,tutteJLMS36} \label{thm:tutte} 
A graph $G=(V,E)$ has $k$ edge-disjoint spanning trees if and only if, for every 
partition ${\mathfrak F}$ of $V$, $|E_G({\mathfrak F})|\geq k(|{\mathfrak F}|-1)$.
\end{theorem}

Let $D = (V,A)$ be a digraph and $r$ be a vertex of $D$.  An {\bf
  out-branching} (respectively, {\bf in-branching)} in $D$ is a
spanning subdigraph, denoted $B^+_r$ (respectively, $B^-_r$), of $D$ in which
each vertex $v \neq r$ has precisely one entering (respectively,
leaving) arc and $r$ has no entering (respectively, leaving) arc. The
vertex $r$ is called the {\bf root} of $B^+_r$ (respectively,
$B^-_r$).
Edmonds characterized digraphs with $k$-arc-disjoint out-branchings.
\begin{theorem}\cite{edmonds1973}
    \label{thm:edmbranch}
   Let  $D=(V,A)$ be a digraph and $k$ a natural number. Then $D$ has $k$ arc-disjoint out-branchings rooted at $r$ if and only if
   \begin{equation}
     \label{in-cond}
      d^-(X)\geq k\hspace{3mm}\forall X\subseteq V-r, X \not= \emptyset,
    \end{equation}
    \noindent{}where $d^-(X)$ denotes the number of arcs entering $X$ in $D$.
  \end{theorem}
  \noindent{}Lov\'asz \cite{lovaszJCT21} gave an algorithmic proof of Theorem \ref{thm:edmbranch} which leads to a polynomial algorithms that  either  constructs the desired branchings or finds a subset violating (\ref{in-cond}).

Thomassen proved that the problem of deciding
whether a digraph contains an out-branching and and
in-branching which are arc-disjoint is NP-complete (see \cite{bangJCT51}). It was proved in \cite{bangJGT42} that this problem is
polynomial time solvable for acyclic digraphs. Furthermore,  acyclic digraphs which contain a pair of arc-disjoint branchings $B^+_s,B^-_t$
 admit a nice characterization \cite{bangJGT42}.

\begin{theorem}\cite{bangJGT42}
  \label{thm:acyclicin-outbr}
  There exists a polynomial algorithm $\cal  A$ which given an acyclic digraph $D$ and two vertices $s,t$ of $D$; decides whether $D$ has arc-disjoint branchings $B^+_s,B^-_t$ and outputs such a pair when they exist and otherwise outputs a certificate showing that there is no such pair in $D$.
  \end{theorem}
 We shall not need the certificate mentioned above but only the following facts.
Suppose that $D = (V,A)$ is an acyclic digraph with an acyclic ordering $v_1,v_2,\ldots{},v_n$ of its vertices
and that $B^+_s,B^-_t$ are 
arc-disjoint branchings rooted at $s,t$ respectively 
in $D$. Then $s$ must be  $v_1$ and this is the unique vertex of in-degree zero and $t$ must be $v_n$ which is the unique 
vertex of out-degree zero in $D$. Furthermore, for every vertex $v_i$ with $i>1$, the unique arc $v_jv_i$ of $B^+_s$ entering $v_i$ satisfies that  $j<i$. Similarly for every vertex $v_i$ with $i<n$ the unique arc $v_iv_q$ leaving $v_i$ in $B^-_t$ satisfies that $q>i$.

Every graph has an acyclic orientation. A natural way of obtaining an acyclic
orientation of a graph $G$ is to orient the edges according to total order 
$\leq $ of $V(G)$, that is, each edge $uv$ of $G$ is oriented from $u$ to $v$ if 
and only if $u <v$. In fact, every acyclic orientation of $G$ can be obtained 
in this way.
Given a vertex ordering $\leq $ of $G$, we use $D_{\leq}$ to 
denote the acyclic orientation of $G$ resulting from $\leq $, and call $\leq$  
{\bf good} if $D_{\leq}$ contains an out-branching and 
an in-branching which are arc-disjoint.

Conversely we say that an acyclic orientation $D$ of a graph $G$ is {\bf good} if $D$ contains an out-branching $B^+_s$ and an in-branching $B^-_t$ that are arc-disjoint. By Theorem \ref{thm:acyclicin-outbr} we can check whether a given total order/ acyclic digraph is good in polynomial time and produce the desired branchings if the answer is yes.

 \section{Notation and Preliminaries}
  Notation not introduced here will be consistent with \cite{bang2009}.
All graphs $G$ considered here are supposed to be finite and loopless but may contain multiple edges.
A graph is {\bf simple} if it has no pair of parallel edges.
An {\bf acyclic} ordering of an acyclic digraph $D=(V,A)$ is an ordering $v_1,v_2,\ldots{},v_n$ of its vertices such that $v_iv_j\in A$ implies that $i<j$.

An {\bf $(s,t)$-triple} of $G$ is 
a triple $(\leq,I,O)$ consisting of a total order $\leq$ of $V(G)$ and two edge disjoint spanning trees $I,O$ of $G$
such that every vertex except $t$ has a larger neighbor in $I$ and every vertex except $s$ has a smaller neighbor in $O$.
That is, if we orient the edges from its smaller endpoint to the other, then every vertex but $t$ has an out-neighbor in $I$
and every vertex but $s$ has an in-neighbor in $O$, so that $I$ becomes an in-branching rooted at $t$ and $O$ becomes
an out-branching rooted at $s$.

Recall that a {\bf generic circuit} is a 2T-graph that has no proper subgraph which is also a 2T-graph. Generic circuits, also sometimes called $M$-circuits, play an important role in rigidity theory for graphs. See e.g. \cite{bergJCT88,lamanJEM4}

\begin{lemma}\cite{bangGOpaper}
  \label{lem:GCintersect}
  Let $G_1=(V_1,E_1)$ and $G_2=(V_2,E_2)$ be  distinct generic circuits of a 2T-graph $G$. Then $E_1\cap E_2=\emptyset$ and $|V_1\cap V_2|\leq 1$
\end{lemma}

The following result shows that not only do all generic circuits have an $(s,t)$-triple but we can choose the vertices $s,t$ of an  $(s,t)$-triple freely (as long as they are distinct). 

\begin{theorem}\cite{bangGOpaper}
  \label{thm:genericgood}
  Let $G$ be a generic circuit.
  For any pair of vertices $s \not= t$ and every edge $e$ incident with $s$ or $t$ there exists
  an $(s,t)$-triple $(\leq,I,O)$ with $e \in E(I)$ and
  an $(s,t)$-triple $(\leq,I,O)$ with $e \in E(O)$. Furthermore, there exists a polynomial algorithm $\cal B$ which given a generic circuit $G$, vertices $s\neq t$ and an edge $e$ incident to $s$ or $t$, produces the desired $(s,t)$-triple.
\end{theorem}

\begin{theorem}\cite{bangGOpaper}
  \label{thm:findallGC}
  There exists a polynomial algorithm $\cal  A$ which given a 2T-graph $G=(V,E)$ as input finds the collection $G_1,G_2,\ldots{},G_r$, $r\geq 1$ of generic circuits of $G$.
\end{theorem}

\begin{corollary}\cite{bangGOpaper}
There exists a polynomial algorithm for deciding whether a 2T-graph $G$ is a generic circuit.
\end{corollary}

\begin{theorem}
  \label{decomposeintoGC}\cite{bangGOpaper}
  There exists a polynomial algorithm for deciding whether the vertex set of a 2T-graph $G=(V,E)$ decomposes into vertex disjoint generic circuits. Furthermore, if there is such a decomposition, then it is unique.
\end{theorem}

\begin{theorem}\cite{bangGOpaper}
    \label{thm:matchingcase}
    Let $G=(V,E)$ be a 2T-graph whose vertex set decomposes into vertex disjoint generic circuits $G_1,G_2,\ldots{},G_r$, $r\geq 1$ and let $E'$ be the set of those edges of $G$ which connect different $G_i$'s.
    If $E'$ is a matching, then $G$ has a good
    orientation.
  \end{theorem}

\begin{theorem}\cite{nashwilliamsJLMS39}
\label{thm:NWcover2T}
The edge set of a graph $G$ is the union of two forests if and only if 
\begin{equation}
  \label{sparse}
|E(G[X])|\leq 2|X|-2
\end{equation}
\noindent{}for every non-empty subset $X$ of $V$.
\end{theorem}

\section{Quartics}
Throughout, a {\bf quartic} is a simple 2T-graph in which all vertices have degree $3$ or $4$.
It is immediate that a quartic has exactly four vertices of degree $3$, called its {\bf transits}.
A quartic is {\bf excellent} if it admits an $(s,t)$-triple for every choice of distinct transits $s,t$.
Every generic circuit of a quartic is itself a quartic, induced, and excellent by Theorem \ref{thm:genericgood}.

\begin{figure}[H]
\begin{center}
\scalebox{0.7}{\includegraphics{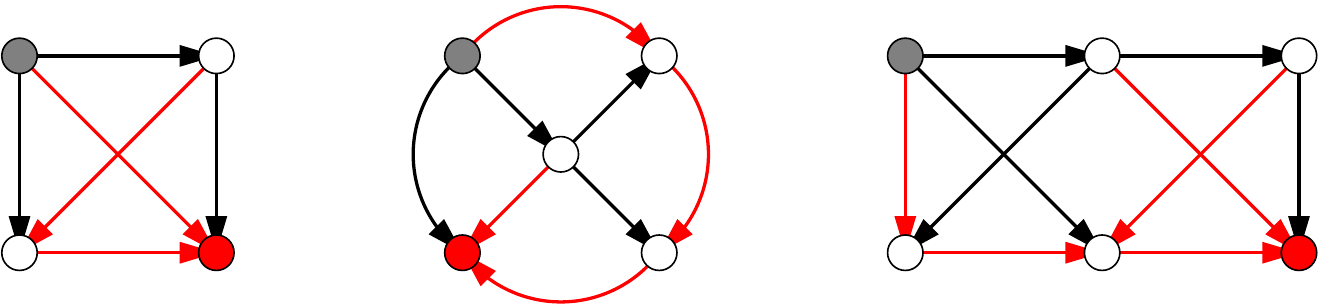}}
\caption{Three quartics, each coming with some $(s,t)$-triple. The source transits $s$ and the $O$'s are displayed black, the sink transits $t$ and the $I$'s are red.
From the indicated acyclic orientation we can derive an ordering $\leq$ of the vertices as required (just take any acyclic ordering of the acyclic digraph).
All three quartics are generic circuits (which ensures the existence of $(s,t)$-triples for every pair of distinct vertices $s,t$).}\label{fig:f1}
\end{center}
\end{figure}

The following is a direct consequence of the definition of a quartic and the fact that these are 2T-graphs.

\begin{proposition}\label{prop:nooftransits}
  Let $Q$ be a proper subquartic of a quartic $G=(V,E)$ and let $d(Q)$ be the number of edges between $V(Q)$ and $V-Q$. Then $d(Q)\in \{2,3,4\}$. Furthermore, the number of transits of $G$ that are vertices of $Q$ is precisely $4-d(Q)$.
  \end{proposition}

It follows from Lemma \ref{lem:GCintersect} that any pair of distinct generic circuits of a quartic must be vertex disjoint. Thus the following is a direct consequence of Theorem \ref{thm:findallGC}
\begin{lemma}
  \label{lem:quarticdecompose}
  Let $G=(V,E)$ be a quartic. Then there is a unique decomposition $V=V_1\cup{}\ldots\cup{}V_r\cup{}V_{r+1}\cup\ldots\cup{}V_p$ such that $r\geq 1$, $G[V_i]$ is a generic circuit with at least 4 vertices for all $i\in [r]$ and each $V_j$ is a single vertex for $j=r+1,\ldots{},p$. Furthermore, the decomposition above can be found in polynomial time.
  \end{lemma}

Suppose that $Q$ and $R$ are disjoint quartics and let $a \not= b$ and $c \not= d$ be two transits of $Q$ and $R$, respectively.
Then the graph $G$ obtained from $Q \cup R$ by adding an edge connecting $a,c$ and an edge connecting $b,d$
is called a {\bf sum} of $Q$ and $R$. Let $X$ be the set of transits of $Q$ or $R$ distinct from $a,b,c,d$.
It is immediate that $G$ is a simple 2T-graph where the four vertices from $X$ have degree $3$ and all others have degree $4$,
in other words: a quartic.

\begin{lemma}
  \label{lem:L1}
  Any sum of two excellent quartics is excellent. 
\end{lemma}

{\bf Proof.}
Let $Q,R,a,b,c,d,G,X$ be just as above and take $s \not= t$ from $X$. By symmetry we may assume that $s \in V(Q)$.

\begin{figure}[H]
\begin{center}
\scalebox{0.6}{\includegraphics{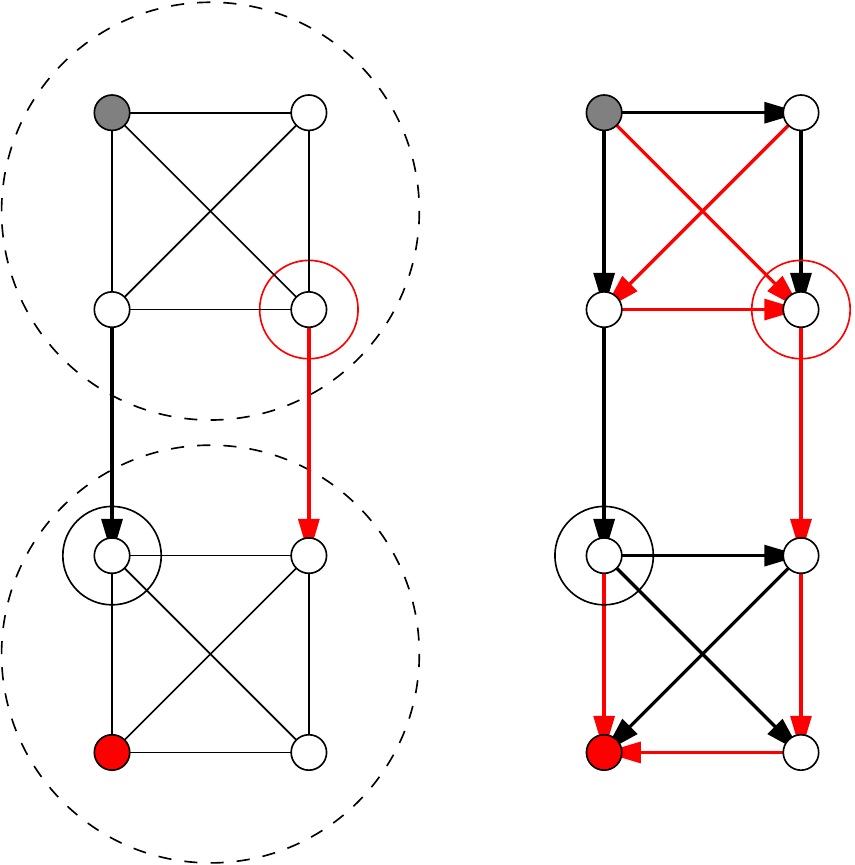}}
\caption{The case that $s,t$ in the proof of Lemma \ref{lem:L1} belong to different summands; $s$ is displayed in black, $t$ in red.
One arbitrarily promises the two connecting edges to $O$ (to be drawn black) and $I$ (in red),
respectively, and think of them as oriented from the source summand to the sink summand.
This determines a local sink $t'$ in the summand containing $s$ (encircled red), and a local source $s'$ in the other part, as in the left side of the picture.
We then find an $(s,t')$-triple and a $(s',t)$-triple  for the summands accordingly and combine as to get what is displayed right.}\label{fig:ff2}
\end{center}
\end{figure}

If $t \in V(R)$ then we take an $(s,b)$-triple $(\leq,I,O)$ for $Q$ and a  $(c,t)$-triple $(\leq',I',O')$ for $R$,
and place the vertices of $V(R)$ right after $b$ in $\leq$ (respecting $\leq'$). Together with $I+I'+bd$ and $O+O'+ac$,
we get an $(s,t)$-triple for $G$. See Figure \ref{fig:ff2}.

\begin{figure}[H]
\begin{center}
\scalebox{0.6}{\includegraphics{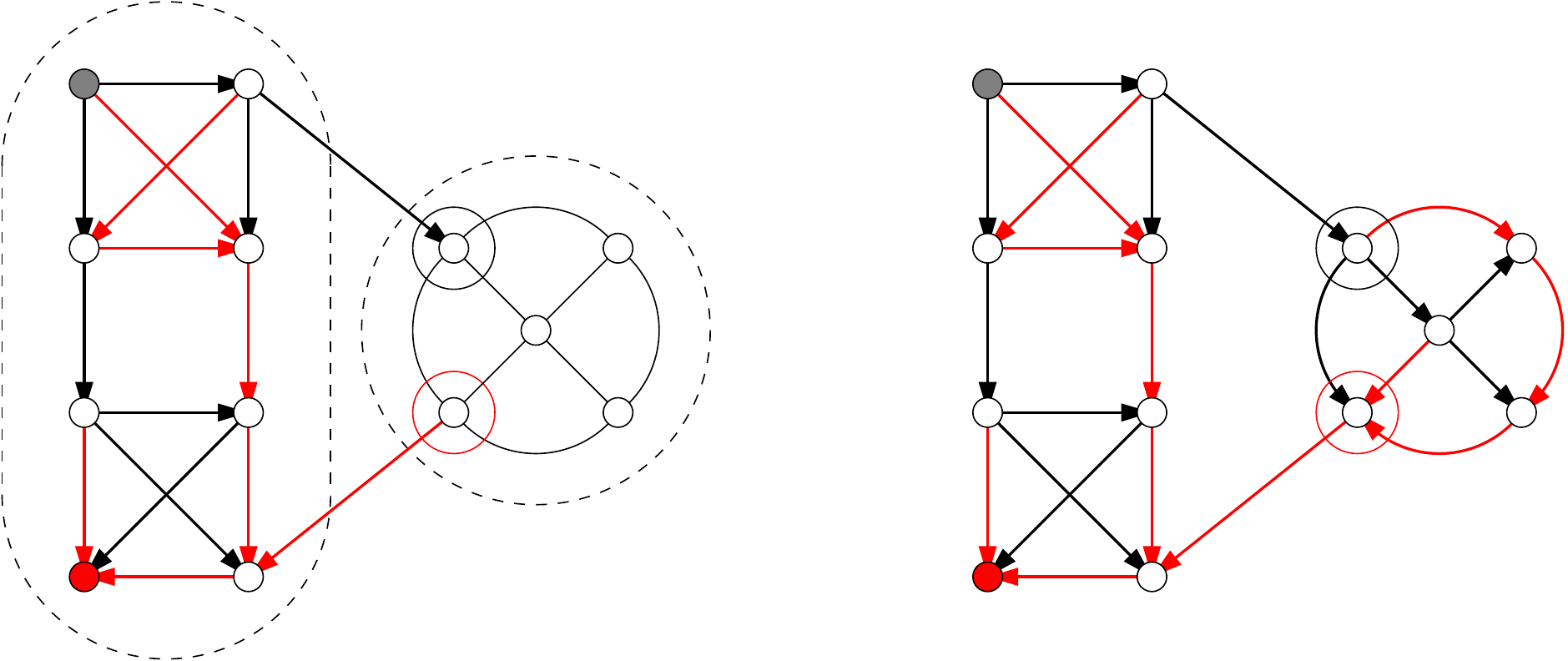}}
\caption{The case that $s,t$ in the proof of Lemma \ref{lem:L1} belong to the same summand;
  $s$ is displayed in black, $t$ in red.
One first calculates an $(s,t)$-triple for that summand (incidentally, it is the one we get from the previous figure)
and compares its two vertices incident with the connecting edges. The edge from the smaller one
to the other summand is supposed to be in $O$ (drawn black), the larger one in $I$ (red), and the end vertices of the connecting edges in the other summand determine a local source $s'$
and a local sink $t'$ there (encircled black and red, respectively, all in the left part of the picture). Based on that, one calculates an $(s',t')$-triple
for the second summand and combines as to get what is displayed right.}\label{fig:f3}
\end{center}
\end{figure}

Next look at the case that $t \in V(Q)$, too.
Since $Q$ is excellent, there exists an $(s,t)$-triple $(\leq,I,O)$ for $Q$.
Without loss of generality, $a<b$ (otherwise we change the roles of $a,b$ and of $c,d$).
Now $R$ has a $(c,d)$-triple $(\leq',I',O')$, and by inserting the vertices of $V(R)$ into $\leq$ in between $a$ and $b$ (respecting $\leq'$)
we get, together with $I+I'+bd$ and $O+O'+ac$, an  $(s,t)$-triple for $G$. See Figure \ref{fig:f3}.\hspace*{\fill}$\Box$

Suppose our graphs are modelled as triples $G=(V,E,\gamma)$, where $V,E$ are just sets and $\gamma \subseteq V \times E$ is
the incidence relation, that is, for every $e \in E$ there exist exactly two vertices $x,y\in V$ such that $x \gamma e$ and $y \gamma e$.
For a partition $\mathfrak{P}$ of $V(G)$, let $E_G(\mathfrak{P})$ denote the set of edges of $G$ with endpoints in two distinct members of $\mathfrak{P}$.
The {\bf quotient graph} $G/\mathfrak{P}$ is then defined to be $(\mathfrak{P},E_G(\mathfrak{P}),\gamma')$,
where $Q \gamma' e$ holds for $Q \in \mathfrak{P}$ if and only if $e$ has an endpoint in $Q$.
The reason why we go down to this level of detail is that we want to ensure that the edges of $G/\mathfrak{P}$ are
actually edges of $G$ (not just by means of correspondence).

\begin{definition}
  \label{def:goodquartic}
Let $G$ be a quartic. 
A partition $\mathfrak{P}$ of $V(G)$ is {\bf excellent} if every set in $\mathfrak{P}$ that contains  more than one vertex induces an excellent quartic.
\end{definition}

\begin{lemma}
  \label{lem:L2}
  Let $G$ be a quartic and let $\mathfrak{P}$ be an excellent partition of $V(G)$. If $G/\mathfrak{P}$ is an excellent quartic then so is $G$.
\end{lemma}

{\bf Proof.} Suppose that $Q$ is an excellent subquartic of $G$ induced by some member of $\mathfrak{P}$.
By definition of a quartic, every edge incident with the vertex $V(Q)$ in $G/\mathfrak{P}$ connects a transit from $Q$,
so Proposition \ref{prop:nooftransits} implies that $V(Q)$ contains a unique transit of $G$ if $V(Q)$ is a transit of $G/\mathfrak{P}$ and no transit of $G$ otherwise.
It follows that for distinct transits $s,t$ of $G$ there exist distinct sets $S,T$ in $\mathfrak{P}$ with $s \in S$ and $t \in T$.
Since $G/\mathfrak{P}$ is excellent, it admits an $(S,T)$-triple $(\leq,I,O)$.
For every vertex $X \in G/\mathfrak{P}$ distinct from $S$ there exists a unique edge $e$ from $O$ connecting it to a smaller vertex,
and we denote its end vertex (as an edge in $G$) in $X$ by $s_X$.
Likewise, for every vertex $X \in G/\mathfrak{P}$ distinct from $T$ there exists a unique edge $f$ from $I$ connecting it to a larger vertex,
and we denote its end vertex (as an edge in $G$) in $X$ by $t_X$. 
Additionally, we set $s_S:=s$ and $t_T:=t$. Since  $G[X]$ has minimum degree 3 and for each $|X|>1$, $s_X$ and $t_X$ are distinct transits of $G[X]$, so, by definition of an excellent quartic,
$G[X]$ has an $(s_X,t_X)$-triple $(\leq_X,I_X,O_X)$. When  $|X|=1$, we let $\leq_X$ be the unique linear order on $X$ and
$I_X,O_X$ be the edgeless graph on $X$. By replacing $X$ in $\leq$ with its members (respecting $\leq_X$) we get a
linear order of $V(G)$, and together with $\bigcup I_X + E(I)$ and $\bigcup O_X+E(O)$ (unions taken over $X \in \mathfrak{P}$)
we get an $(s,t)$-triple for $G$.
\hspace*{\fill}$\Box$

\begin{figure}[H]
\begin{center}
\scalebox{0.6}{\includegraphics{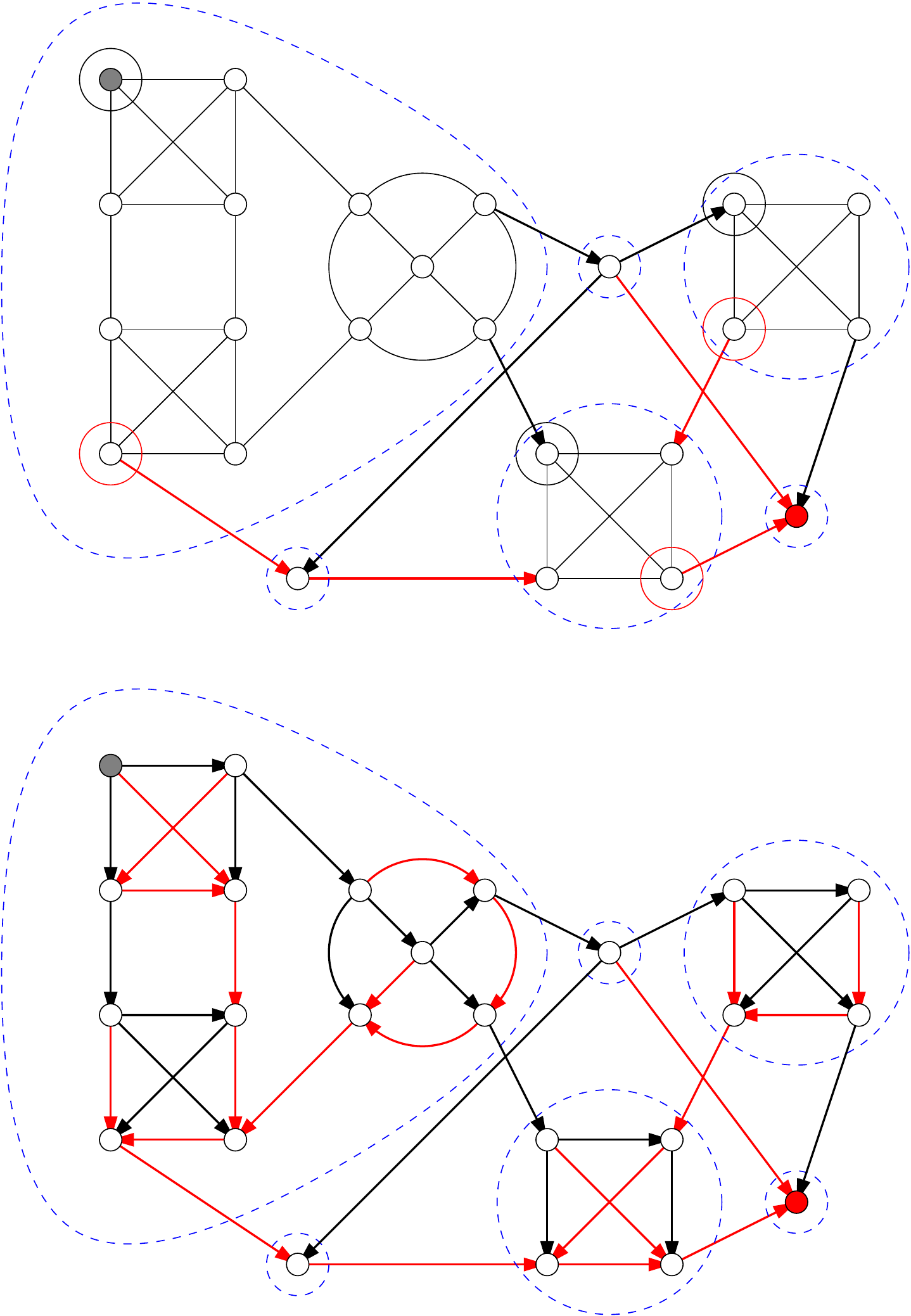}}
\caption{The picture illustrates the proof of Lemma \ref{lem:L2} and shows an excellent partition of some quartic into six pieces in the upper part;
the source transit is displayed black, the sink transit is red.
In this example the quotient quartic is a generic circuit and, thus, by Theorem \ref{thm:genericgood} it is  excellent; for the two classes containing the source and sink, respectively, one first calculates
a triple of the quotient graph (we did it here according to the rightmost part of Figure \ref{fig:f3}).
This determines local sources and sinks in the partition classes (encircled black and red for the nontrivial classes).
We then calculate triples for these accordingly and combine as
to obtain what is displayed in the lower part.}\label{fig:f4}
\end{center}
\end{figure}

We call a quartic {\bf normal} if the edge neighborhood of every proper subquartic is a matching of size $3$ or $4$.
By Theorem \ref{thm:NWcover2T} every 2T-subgraph of a 2T-graph is an induced subgraph. Hence if $Q$ is a subquartic of $G$ then every subquartic $Q'$ of $Q$ is again a subquartic of $G$. 

\begin{theorem}
  \label{thm:normalisgood}
  Every normal quartic is excellent.
\end{theorem}

{\bf Proof.}
Let $G$ be a normal quartic.
The partition of $V(G)$ into singletons is an excellent partition, so  there exists a coarsest excellent partition, say, $\mathfrak{P}$.
It suffices to prove that $\mathfrak{P}=\{V(G)\}$, so suppose that this does not hold.

It is easy to check that $G/\mathfrak{P}$ is a 2T-graph, and since $G$ is normal, all vertices have degree $3$ or $4$.
If $G/\mathfrak{P}$ contains a pair $e,f$ of parallel edges then, as $G$ is normal, $e,f$ form a matching in $G$ connecting subquartics
of $G$ induced by distinct  members $X,Y$ of $\mathfrak{P}$. Note that the sum $G[X\cup Y]$ of $G[X]$ and $G[Y]$ is a  subquartic of $G$, so by Lemma \ref{lem:L1} it is excellent. In this case we set $\mathfrak{P}^+:=(\mathfrak{P} \setminus \{X,Y\}) \cup \{X \cup Y\}$
and note that  $\mathfrak{P}^+$ is excellent and coarser than $\mathfrak{P}$ .

Otherwise, $G/\mathfrak{P}$ is a quartic. By Lemma \ref{lem:quarticdecompose},  $G/\mathfrak{P}$ has at least one generic circuit. Thus we can take a generic circuit of $G/\mathfrak{P}$, with vertex set $\mathfrak{C} \subseteq \mathfrak{P}$,
and deduce from Theorem \ref{thm:genericgood} that it is an excellent subquartic of $G/\mathfrak{P}$.
In this case we set $\mathfrak{P}^+:=(\mathfrak{P} \setminus \mathfrak{C}) \cup \{\bigcup \mathfrak{C}\}$ and note that
$G[\bigcup \mathfrak{C}]$ is a subquartic of $G$. By Lemma \ref{lem:L2} we know that it is excellent.

In either case, $\mathfrak{P}^+$ is an excellent partition coarser than $\mathfrak{P}$, contradiction.
\hspace*{\fill}$\Box$

\begin{figure}[H]
\begin{center}
\scalebox{0.6}{\includegraphics{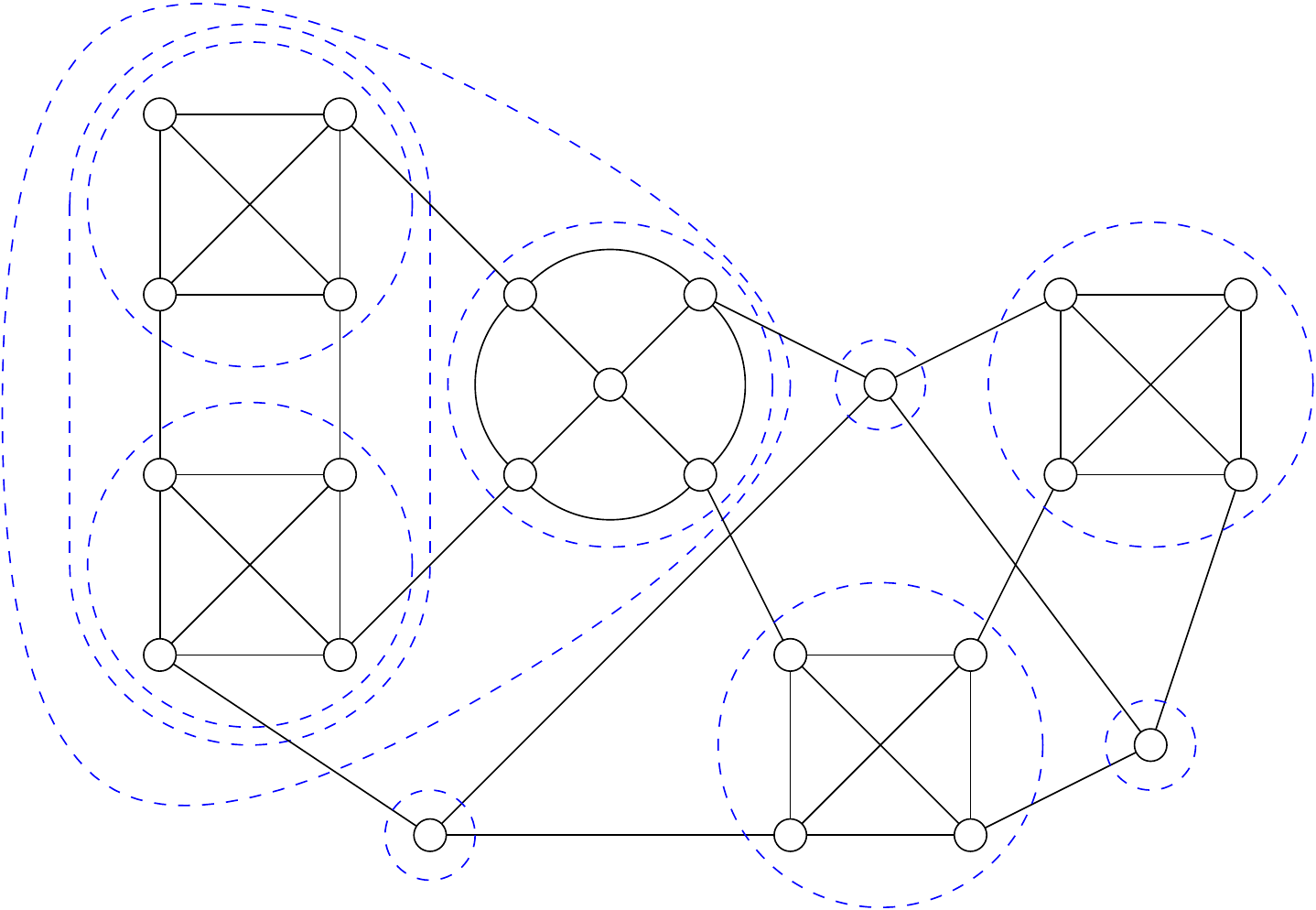}}
\caption{How to get coarser and coarser excellent partitions as in the proof of Theorem \ref{thm:normalisgood}.
By Lemma \ref{lem:quarticdecompose} we can immediately obtain a partition which is coarser than the partition into singletons by partitioning into singletons and generic circuits (the latter are known to be disjoint and there is at least one);
in the example, this is the partition given by the eight dotted circles, and we recognize that the classes formed by the two leftmost $K_4$'s are connected by
two parallel edges in the quotient graph, so that we can coarsen by identifying these two classes. In the new quotient graph,
the new class is connected by two parallel edges to the class to its right formed by a $4$-wheel,
so that we can coarsen again. We get the partition of Figure \ref{fig:f4}, and so Lemma \ref{lem:L2}
proves that $\{V(G)\}$ is an excellent partition.}\label{fig:f5c}
\end{center}
\end{figure}

We will now show that normal quartics can be recognized in polynomial time. To do so we need some results on intersections of quartics.

\begin{lemma}
  \label{lem:normalis3ec}
  Every normal quartic is 3-edge-connected.
\end{lemma}
\pf Suppose that $G=(V,E)$ is a quartic that is not 3-edge-connected and let $X,V-X$ be a partition such that $d(X)=2$ and every nonempty proper subset of $X$ has at least 3 neighbours.
Then  $|X|>1$ as $G$ is a quartic and $G[X]$ is a 2T-graph and it is also a quartic by the minimality of $X$. Now $X$ shows that $G$ is not normal.
\qed

\vspace{2mm}

To simplify notation, when we speak of a quartic $Q$  below we will sometimes think of $Q$ as an induced subgraph and sometimes as a set of vertices.
We call two quartics $Q,Q'$ {\bf skew} if their intersection $Q\cap Q'$ is a non-empty proper subset of both. 

\begin{lemma}
  \label{lem:skewQ}
  Let $L$ and $R$ be skew subquartics of a normal quartic $G$. Then the following holds.
  \begin{itemize}
    \item[(a)] $L\cup R$ is a quartic.
  \item[(b)] $G$ has exactly 2 edges $aa',bb'$ with $a,b\in L\cap R$, $a',b'\in L-R$ and  exactly two edges
    $cc',dd'$ with $c,d\in L\cap R$ and $c',d'\in R-L$.
    \item[(c)] There is no edge in $G$ from $L-R$ to $R-L$.
    \item[(d)] $L-R,R-L,L\cap R$ are all 2T-graphs.

    \end{itemize}
  \end{lemma}
  \pf As $R,L$ are 2T-graphs there are at least two edges from $L\cap R$ to each of $L-R$ and $R-L$.
  Denote the endvertices of those in $L\cap R$ by $a,b$ and $c,d$, respectively. As $R$ and $L$ are quartics and hence have minimum degree 3, we have $a\neq b$ and $c\neq d$.  Note that  $a,b$ are transits of $R$ and $c,d$ are transits of $L$ .
  We also have $d(x)=4$ for $x\in \{a,b,c,d\}$, implying that $R\cup L$ has at most 4 vertices of degree 3. It also follows from Theorem \ref{thm:NWcover2T} that $R\cup L$ has at least 4 vertices of degree 3, so it has exactly 4 vertices of degree 3 and hence it is a quartic, showing that  (a) holds. To see that (b) holds, suppose that there is another edge $ef$ from $L\cap R$ to $L-R$, then $e$ is a transit of $L$ but not of $L\cup R$, contradicting that $L\cup R$ has 4 transits. This proves the first part of (b) and the second follows analogously. Now we see that (c) holds since if there was an edge $lr$ from $L-R$ to $R-L$, then $l$ ($r$) is a transit of $L$ ($R$) but has degree 4 in $G$ so $L\cup R$ would have at most 2 transits, contradicting Theorem \ref{thm:NWcover2T}. Finally (a), (b) and (c) imply that (d) holds. \qed

  \2
  
  Since generic circuits contain no proper 2T-graphs we get the following consequence of Lemma \ref{lem:skewQ}.
  
  \begin{corollary}
    \label{cor:GCnotskew}
    If $Q,Q'$ are subquartics of a normal quartic $G$ and $Q$ is a generic circuit, then $Q$ and $Q'$ are not skew.
  \end{corollary}

  Let us call a proper subquartic of a quartic {\bf bad} if its edge neighborhood is not a matching of size $3$ or $4$. Obviously every non-normal quartic contains a minimal bad subquartic.

  \begin{lemma}
    \label{lem:minbadskew}
    Let $Q$ be  a minimal bad subquartic of a quartic $G$ such that $a,b\in Q$ have a common neighbour $c$ in $V-Q$. If $Q$ is skew with another subquartic $S$, then  $c\not\in S$ and $S$ contains exactly one of the vertices $a,b$.
  \end{lemma}

  \pf Suppose that $S$ is a subquartic which is skew with $Q$. Suppose first that $a,b\in Q-S$. By Lemma \ref{lem:skewQ}(c) $c\not\in S$ and by (d) $Q-S$ is a 2T-graph. If $Q-S$ is also a quartic, then it is bad (the vertex $c$ certifies this), contradicting the minimality of $Q$. Hence $Q-S$ contains a vertex of degree 2. Now remove vertices of degree 2 from $Q-S$ as long as possible and note that, as every 2T-graph contains a generic circuit,  this will end in a non-empty 2T-graph $Q'$. Let $x$ be the last vertex that we removed to obtain $Q'$ and observe that $x$ certifices that $Q'$ is bad, contradicting the minimality of $Q$. If $a,b$ are both  in $S$ we reach a similar contradiction by finding a smaller bad quartic inside $Q\cap S$. Finally, by Lemma \ref{lem:skewQ} (c), we cannot have $S\cap \{a,b,c\}=\{a,c\}$ or $S\cap \{a,b,c\}=\{b,c\}$. \qed

  \2

Now we are ready to prove that normal quartics can be recognized in polynomial time.

\begin{theorem}
  \label{thm:decidenormal}
  There exists a polynomial algorithm for deciding whether a given quartic is normal. Furthermore, if the input is not normal, then the algorithm will produce a certificate for this.
  \end{theorem}

  \pf Let us first see that if the input $G$ is a normal quartic,  then we can make  a polynomial coarsification algorithm ${\cal A}_C$  by following the steps of the proof of Theorem \ref{thm:normalisgood} as follows. Start by applying algorithm ${\cal A}'$ from Lemma \ref{lem:quarticdecompose} to obtain the partition $\mathfrak{P}_0= \{Q_1,\ldots{},Q_r,\{v_{r+1}\},\ldots{},\{v_p\}\}$ of $G$ into generic circuits and singletons and form the quotient $G/\mathfrak{P}_0$.
  Denote by $G/\mathfrak{P}_i$ the quotient after the $i$th coarsification step.
  If $G/\mathfrak{P}_i$ has a pair of parallel edges, then they must go between distinct members $X,Y$  of $\mathfrak{P}_i$ none of which are singletons (as $G$ is normal) and we obtain $\mathfrak{P}_{i+1}$ by replacing  $X,Y$  by $X\cup Y$ to obtain a coarser partition. Otherwise we apply the algorithm of Lemma \ref{lem:quarticdecompose} to
  the quartic  $G/\mathfrak{P}_i$ to find a generic circuit $Z$ of this graph and then obtain  $\mathfrak{P}_{i+1}$ by taking the union of those subsets $\mathfrak{P}_i$ which correspond to $Z$ in $G/\mathfrak{P}_i$.
  This continues until the new partition has only one set.\\

  Now let $G$ be a given quartic.
  If $G$ is  only 2-edge-connected (has a 2-edge cut), then it is bad by Lemma \ref{lem:normalis3ec} and we can easily detect this in polynomial time so we may assume that the input $G$ is 3-edge-connected quartic.
  Suppose that $G$ is not normal and let $Q$ be a minimal bad subquartic. We  claim that some bad quartic will appear as one of the subsets of one of the $\mathfrak{P}_i$'s so that, by checking whether the edge neighbourhood of each part of the current partition is a matching of size 3 or 4,  ${\cal A}_c$ can  detect $Q$ or another bad quartic and use this  as a certificate that $G$ is not normal.\\

  As $Q$ is a bad quartic it must contain two distinct vertices $a,b$ with a common neighbour $c$ in $V-Q$ (certifying that it is bad).  Consider first the initial
  partition $\mathfrak{P}_0$. By Corollary \ref{cor:GCnotskew} we cannot have that $a,b,c$ all belong to one of the generic circuits $Q_1,\ldots{},Q_r$. If some $Q_i$ contains $a,b$ but not $c$, then $Q_i$ certifies that $G$ is bad (possibly $Q_i=Q$). Now it follows from  Corollary \ref{cor:GCnotskew}   that  if one of $Q_1,\ldots{},Q_r$ intersects $\{a,b,c\}$ then it contains precisely one of $a,b,c$. We prove our claim by tracking what happens to the subsets $A,B,C$ of the current partition $\mathfrak{P}_i$ which contain $a,b,c$, respectively, as the coarsification progresses. By the previous remarks, after the first partition $\mathfrak{P}_0$ is created, the sets $A,B,C$ are still disjoint. Let $i>0$ be the step in which at least two of the sets $A,B,C$ are merged  when we   move from $\mathfrak{P}_{i-1}$ to $\mathfrak{P}_i$. 
  It follows from Lemma \ref{lem:minbadskew} and the fact that all non-trivial sets in each $\mathfrak{P}_i$ are
  quartics that, in fact,  all of $A,B,C$ are merged when we move from $\mathfrak{P}_{i-1}$ to $\mathfrak{P}_i$.  First observe that at this point $C=\{c\}$ since if $\{c\}$ would be part of a generic circuit $Z$ of $G/\mathfrak{P}_j$ for some $j<i$, then the union of the sets of $\mathfrak{P}_j$ corresponding to $Z$ would have been detected earlier as a bad quartic.\\

  Let  $W$ be the  generic circuit of $G/\mathfrak{P}_{i-1}$ that the algorithm uses to replace $\mathfrak{P}_{i-1}$ by $\mathfrak{P}_i$, let $\mathfrak{M}$ be the subpartition of $\mathfrak{P}_{i-1}$ whose sets correspond to the  generic circuit $W$ of $G/\mathfrak{P}_{i-1}$ and let  $Q_W$ be the quartic whose vertex set is the union of all the sets in $\mathfrak{M}$. By Lemma \ref{lem:skewQ}, we have $Q\subseteq Q_W$.
  Suppose first that $A,B$ are both subsets of $Q$ and recall that, as $C=\{c\}$, we have $C\cap Q=\emptyset$.
  It follows from Corollary \ref{cor:GCnotskew} and the fact that $c\not\in Q$ that $Q$ corresponds to a proper subpartition $\mathfrak{N}$ of $\mathfrak{M}$. But now, the fact that $Q$ is a 2T-graph implies that  $G/\mathfrak{N}$ is a 2T-subgraph of $G/\mathfrak{M}$, contradicting that this is a generic circuit of $G/\mathfrak{P}_{i-1}$.\\
  
  Thus we may assume that one of $A,B$, say w.l.o.g. $A$, is  skew with $Q$. By Lemma \ref{lem:skewQ}, two of the edges incident to $Q$ in $G$ enter $A-Q$ and $Q$ also has the two edges $ac,bc$ to $C=\{c\}$, hence no other member of $\mathfrak{M}$ is skew to $Q$, in particular $B\subset Q$. Let $u\neq v$ be the two vertices of $Q\cap A$ that send an edge to $A-Q$ and similarly let $x\neq y$ be the two vertices of $Q\cap A$ that send an edge to $Q-A$. So we have identified  5 distinct edges between
  $A\cap Q$ and $V-(A\cap Q)$.  By Lemma \ref{lem:skewQ} $A\cap Q$ is a 2T-graph, so   we get
  $$4|A\cap Q|-4=\sum_{v\in A\cap Q}d_{A\cap Q}(v)\leq 4|A\cap Q|-5,$$
  \noindent{} a contradiction.

  The argument above proves that if there is a bad subquartic in $G$, then some bad subquartic will be detected by ${\cal A}_C$. \qed

  \2
  
  \2

\begin{theorem}
  \label{thm:orientalg}
  There exists a polynomial algorithm $\cal A$ which given a normal quartic $G$ and distinct transits $s,t$ of $G$ returns an $(s,t)$-triple of $G$.
\end{theorem}
\pf We give a high level description which should suffice to show that the algorithm can indeed be implemented as a polynomial algorithm. Whenever we talk about edges, these will correspond to actual edges of $G$ via the relation $\gamma$ that we defined earlier. In order to simplify the description of $\cal A$, we will describe it so that, instead of producing an $(s,t)$-triple of $G$, $\cal A$ will produce an acyclic orientation $D$ of $G$ with  arc-disjoint branchings $B^+_s,B^-_t$. As mentioned in the introduction, this will correspond to an $(s,t)$-triple $(\leq, I, O)$, where $\leq $ is any acyclic ordering of $D$, $I$ is the set of edges that were oriented as arcs in $B^-_t$ and $O$ is the set of edges that were oriented as arcs in $B^+_t$.\\

\begin{figure}[H]
\begin{center}
\scalebox{0.65}{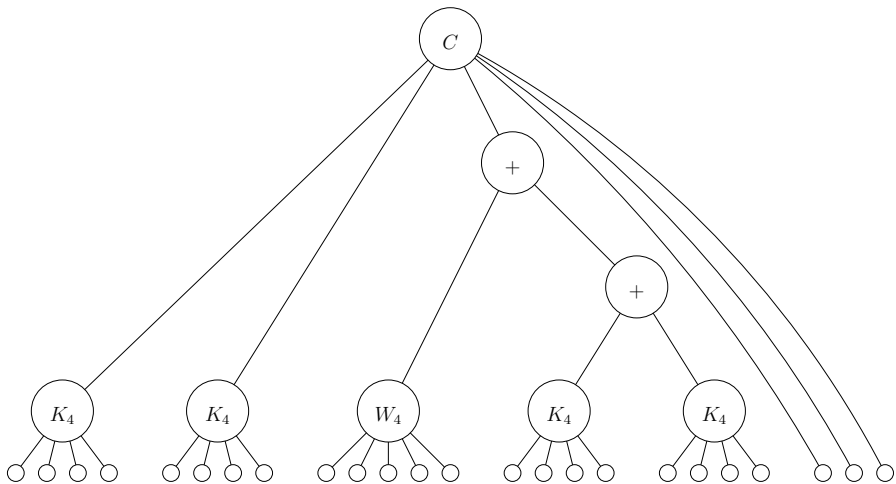}
\caption{A tree $T$ which records a possible result of the coarsification process for the quartic $G$ in Figure \ref{fig:f5c}.
The generic circuit $C$ is the rightmost  graph in Figure \ref{fig:f1} and $W_4$ is the 4-wheel in the middle of Figure \ref{fig:f5c}. There are two kinds of  operators, namely the one indicated by a '+' and the ones which are named by a generic circuit. In the first case, the children of the '+' node are the two subtrees corresponding to the subquartics $X,Y$ which are replaced by their sum. In the second case the children of a vertex labelled by a generic circuit $H$ are those  subtrees which correspond to those singletons and  subquartics which are replaced by one set in the new partition $\mathfrak{P}^+$ after replacing all of these by the union of their vertex sets.}\label{fig:coarsefig}
\end{center}
\end{figure}

Let $G=(V,E)$ be the input quartic and $s,t$ the specified roots. The algorithm $\cal A$  uses the algorithm ${\cal A}'$ of Lemma \ref{lem:quarticdecompose}, starting from the quotient $G/\mathfrak{P}$ formed by the sets of the first decomposition $\mathfrak{P}$ into generic circuits and singletons found by ${\cal A}'$, to execute the coarsification steps outlined in the proof of Theorem \ref{thm:normalisgood}  until it reaches the last step where the current partition $\mathfrak{P}$ satifies that $G/\mathfrak{P}$ is a generic circuit. Note that, as $G$ is normal, the last step in the coarsification process will result from $G/\mathfrak{P}$ being  a generic circuit. This process can be recorded in a tree $T$  as depicted in Figure \ref{fig:coarsefig}. The algorithm $\cal A$ uses the tree $T$ to guide the order in which it produces the orientation as follows. At any time during its execution $\cal A$   will maintain a collection of vertex disjoint trees $T_1,\ldots{},T_k$ with specifications of  local roots $(s_1,t_1),\ldots{},(s_k,t_k)$, respectively, that are currently  unprocessed (we call them {\bf active}), starting with $k=1$ and $T_1=T,s_1=s$ and $t_1=t$ and $\cal A$  will terminate when all trees are processed. $\cal A$ initializes the sets $O$ and $I$ to be empty sets and $A$ to be the empty set ($A$ will contain the edges that are oriented so far). The algorithm uses a stack on the roots of the active trees to control which tree is processed next. So initially the stack contains only one element, namely the root of $T$.

While processing the current tree $T'$ (indicated by the root at the top of the stack) and its local roots $s',t'$, $\cal A$ distinguishes two cases. 

\begin{itemize}
\item The root $R'$ of $T'$ is labelled by a generic circuit $X$. Let $T_{i_1},\ldots{},T_{i_p}$, $p\geq 4$ be the subtrees whose roots are the children of  $R$ in $T'$. Let $T_{i_h},T_{i_q}$ be the trees that contain $s'$ and $t'$ respectively as leaves.
  Note that we cannot have   $i_h=i_q$, that is, $s'$ and $t'$ belong to the same subtree $T_{i_h}$, because then  the subquartic $H$ induced by the leaves of $T_{i_h}$ would be bad by Proposition \ref{prop:nooftransits}.  So $i_h\neq i_q$. Now $\cal A$ picks an arbitrary edge from the edge neighbourhood of the quartic or singleton vertex with corresponds to $T_{i_h}$ and uses the algorithm $\cal B$ of Theorem \ref{thm:genericgood} to find an $(s',t')$-triple $(\leq{}',I',O')$ of $X$. It then orients the edges of $I'$ ($O'$) as an in-tree (out-tree) rooted at $t'$ ($s'$)  and adds these arcs to $I$ ($O$). As illustrated in Figure  \ref{fig:f4},
  this will define local roots in each of the non-trivial subquartics that we indentified to get  $X$ (these 
  correspond to those trees among $T_{i_1},\ldots{},T_{i_p}$ that are not just a leaf of $T'$). We pass on this information by recording together with a non-trivial tree $T_{i_j}$ its two local roots $s_{i_j},t_{i_j}$. Now we mark $T'$ as processed and add each non-trivial tree $T_{i_j}$ and its local roots $s_{i_j},t_{i_j}$ to the active list and add the  roots of  $T_{i_1},\ldots{},T_{i_p}$ to the stack (in any order).
\item The root $R''$ of $T''$ is a '+' node. Let $T_{j_1},T_{j_2}$ be the two children of  $R''$ in $T''$ and let $s'',t''$ be the local roots of $T''$. Let $V_i$ be the sets of leaves of $T_{j_i}$ for $i=1,2$ and let $ac,bd$ denote the two edges in $G$ that go from $V_1$ to $V_2$ (as in the proof of Lemma \ref{lem:L1}). Suppose first  that $s''\in V_i$ and $t''\in V_{3-i}$, w.l.o.g. $i=1$.  Now $\cal A$ orients the edges $ac,bd$  as the arcs $ac$ and $bd$ and adds $ac$ to $O$ and $bd$ to $I$. Then it assigns the local roots $s'',b$ to $T_{j_1}$ and local roots $c,t''$ to $T_{j_2}$, adds these trees to the active list and the stack (in any order) and marks $T''$ as processed.\\
  Suppose next w.l.o.g. that $s'',t''$ are both in $V_1$. This corresponds to the second case in the proof of Lemma \ref{lem:L1}. Here $\cal A$ must process $T_{j_1}$ completely before it can process $T_{j_2}$, so it adds $T_{j_1}$ to the active list and  starts processing $T_{j_1}$. When this terminates, $G[V_{j_1}]$ is oriented as an acyclic digraph $D_{j_1}$ and $\cal A$ now computes an acyclic ordering of this digraph to check whether $a$ is before $b$ in this ordering. If it is, then $\cal A$ orients the edges $ac,bd$ as the arcs $ac,db$, adds $ac$ to $O$, adds $db$ to $I$ and  assigns local roots $s_{j_2}=c,t_{j_2}=d$ to $T_{j_2}$
  Otherwise $a$ is after $d$ in the acyclic ordering and in this case  $\cal A$ orients the edges $ac,bd$ as the arcs $ca,bd$, adds $ca$ to $I$, adds $bd$ to $O$ and  
 assigns local roots $s_{j_2}=d,t_{j_2}=c$ to  $T_{j_2}$. Now $\cal A$ adds $T_{j_2}$ to the active list and the stack and marks $T''$ as processed.
  \end{itemize}

  This completes the description of $\cal A$ and it can easily be checked that it will run in polynomial time. The correctness of $\cal A$ follows from the proof of Theorem \ref{thm:normalisgood}. \qed

  \2
  
Using Lemma \ref{lem:L1} for induction, one can immediately generalize Theorem \ref{thm:normalisgood} to the following result which  generalizes Theorem \ref{thm:matchingcase} in  the case when the 2T-graph $G$ is a quartic.

\begin{corollary}
  \label{C1}
  Let $G$ be a quartic such that the edge neighborhood of every proper subquartic is a matching. Then $G$ is excellent.
\end{corollary}
\pf If $G$ is a generic circuit, then the claim follows from Theorem \ref{thm:genericgood} so we may assume that $G$ is not a generic circuit and that the claim holds for all quartics with fewer vertices than $G$.
Suppose that $G$ contains a proper subquartic $Q$ whose edge neighbourhood is a matching of size 2. Then $Q'=G[V-V(Q)]$ is also a quartic and by induction $Q$ and $Q'$ are excellent. Now it follows from Lemma \ref{lem:L1} that also $G$ is excellent.
Thus we may assume that the edge neighbourhoood of every proper subquartic is a matching of size 3 or 4 (it cannot be larger as $G$ has maximum degree 4) and the claim follows from Theorem \ref{thm:normalisgood}. \qed

\section{Every 4-regular 4-connected graphs is  good}
The following result, which implies that $4$-connected line
graphs of cubic graphs  have a spanning generic circuit, was proved in \cite{bangGOpaper}

\begin{theorem}\cite{bangGOpaper}
  \label{thm:4reg4con}
  Let $G$ be a $4$-regular $4$-connected graph in which every edge is
  on a triangle.  Then $G-\{e,f\}$ is a spanning generic circuit for
  any two disjoint edges $e,f$.  In particular, $G$ admits a good
  ordering.
\end{theorem}

Using the results on quartics above, we can generalize this as follows.

\begin{theorem}
  \label{thm:4reg4confull}
  Let $s \not= t$ be vertices of a $4$-regular, $4$-connected graph $G$.
  Then there exists an $(s,t)$-triple for $G$.
\end{theorem}

{\bf Proof.}
Observe that $G$ must be  simple since it is 4-regular and
4-connected.
There exists a pair of disjoint edges $e,f$ such that $e$ is incident with $s$ and $f$ is incident with $t$. By Theorem \ref{thm:tutte} and the fact that $G$ is 4-connected, $G':=G-\{e,f\}$ is a quartic, and
$V(e) \cup V(f)$ is the set of its transits. 
Let $H$ be a proper subquartic of $G'$ and let $X=V(G') \setminus V(H)$.
As $H$ is a quartic, the edge neighborhood $N'$ of $X$ in $G'$ consists of at most $4$ edges, and their endpoints in $V(H)$ are pairwise distinct transits of $H$ and clearly not transits of $G'$.
It follows that $V(H)$ contains exactly $4-|N'|$ transits of $G'$, so that the edge neighborhood $N$ of $V(H)$ in $G$ consists of 
at most $|N'|+4-|N'|=4$ edges. Since $G$ is $4$-connected, $|N|=4$. As $X$ contains $|N'|$ transits of $G'$, $|X| \geq 2$ holds,
and from regularity and simplicity of $G$ we infer that $|X| \geq 4$ holds. By $4$-connectivity, $V(H)$ must have at least four neighbors in $G$,
so that $N$ and hence $N' \subseteq N$ are matchings. Therefore, Corollary \ref{C1} applies.
\hspace*{\fill}$\Box$

\section{Dense graphs are (super) good}

Let us say that a graph $G$ has property (*), if for all $s \not =t$ from $V(G)$ there exists a spanning 2T-graph admitting an $(s,t)$-triple.
The following is  straightforward to check.

\begin{proposition}
  \label{prop:addone}
  If $G$ has property (*) then so has any graph obtained from $G$ by adding a new vertex $x$
  and edges from $x$ to at least two vertices from $V(G)$.
\end{proposition}

Moreover, if $G,H$ are disjoint graphs having property (*) and $a,b,c,d$
are pairwise distinct vertices with $a,b \in V(G)$ and $c,d \in V(H)$ then the graph obtained from the union of $G,H$
by adding two edges connecting $a,c$ and $b,d$, respectively, will have property (*); the proof is almost literally the proof of Lemma \ref{lem:L1}.

We now look at dense graphs. Let $n \geq 2$.
The union of disjoint copies of $K_{\lfloor n/2 \rfloor}$ and $K_{\lceil n/2 \rceil}$ is a disconnected graph on $n$ vertices with minimum degree $\lfloor n/2 \rfloor-1$,
showing that a minimum degree of at least $\lfloor n/2 \rfloor$ is necessary as to guarantee that a graph is connected; this bound is also
sufficient.

In fact a simple graph $G$ on $n$ vertices with $\delta(G) \geq \lfloor n/2 \rfloor$ will be already $\lfloor n/2 \rfloor$-edge connected:
To see this, take any $X \subseteq V(G)$ with $1 \leq |X| \leq \lfloor n/2 \rfloor$ and observe that there 
are at least $|X| \cdot (\lfloor n/2 \rfloor - |X|+1)$ edges connecting some vertex from $X$ to some vertex from $V(G) \setminus X$;
this bound takes its minimum, $\lfloor n/2 \rfloor$, for $|X|=1$ and $|X|=\lfloor n/2 \rfloor$.
For $n \geq 8$, every  graph $G$ on $n$ vertices with $\delta(G) \geq \lfloor n/2 \rfloor$  is $4$-edge-connected and, thus, by Theorem \ref{thm:tutte}, contains a spanning 2T-graph, which is clearly necessary for having property (*).
For $n<8$, there exist simple graphs $G$ on $n$ vertices $\delta(G) \geq \lfloor n/2 \rfloor$ which are simply too sparse
to admit a spanning 2T-graph. However, if they do, then $n \geq 4$ and we show that they have property (*), with one exception.

\begin{lemma}
  \label{lem:L1n}
  Suppose that $n \leq 7$ and $G$ is a simple graph on $n$ vertices with $\delta(G) \geq \lfloor n/2 \rfloor$.
  If $G$ has a spanning 2T-graph then it has property (*) unless $G$ is the graph on
  seven vertices obtained by identifying two disjoint copies of $K_4$ at a single vertex.
\end{lemma}

{\bf Proof.} Let $C$ be a generic circuit of $G$; as $C$ has property (*) by Theorem \ref{thm:genericgood}, it extends to a maximal subgraph $H$ of $G$ having property (*).
Clearly, $n \geq |V(H)| \geq |V(C)| \geq 4$. If $H$ spans $G$ then we are done, otherwise $X:=V(G) \setminus V(H)$ is not empty.
By Proposition \ref{prop:addone} and the  maximality of $H$, (**) every vertex from $X$ has at most one neighbor in $V(H)$.
Since $G$ contains a spanning 2T-graph, $\delta(G) \geq 2$ so  $|X| \geq 2$, and if $|X|=2$ then $n \geq 6$ and $\delta(G) \geq 3$, implying that every vertex in $X$ has
two neighbors in $V(H)$, contradiction. So $|X|=3$, $n=7$, $H=C \cong K_4$, and (**) implies that $G[X]$ is a triangle.

\begin{figure}[htbp]
\begin{center}
\scalebox{0.6}{\includegraphics{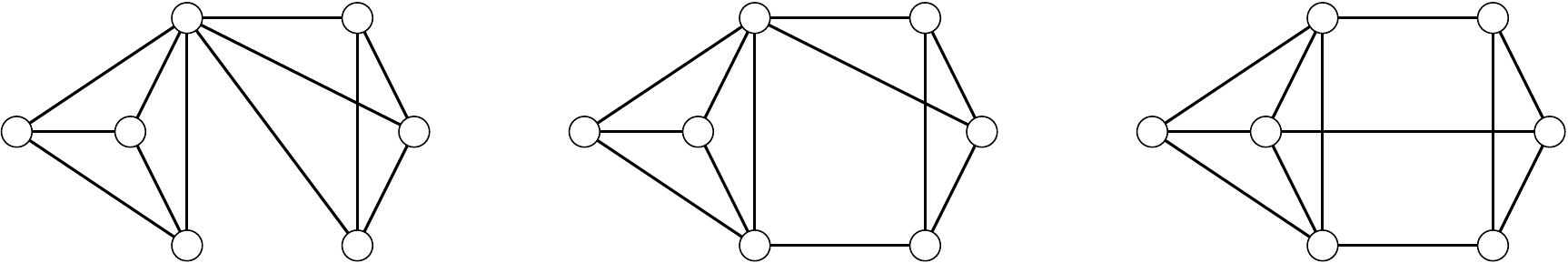}}
\caption{The three remaining graphs to investigate. The leftmost is the exceptional graph, the others have property  (*).}\label{fig:f3n}
\end{center}
\end{figure}

Since $\delta{}(G)\geq 3$ the maximality of $H$ implies that
every vertex in $X$ has exactly one neighbor in $V(H)$. Depending on the size of the neighborhood of $X$ in $G$ ($1$, $2$, or $3$),
we get one of the three graphs in Figure \ref{fig:f3n}. The leftmost is the exceptional graph, so let us concentrate on the other two.
We will prove that they have property (*), contradicting the maximality of $H$. 

\begin{figure}[H]
\begin{center}
\scalebox{0.6}{\includegraphics{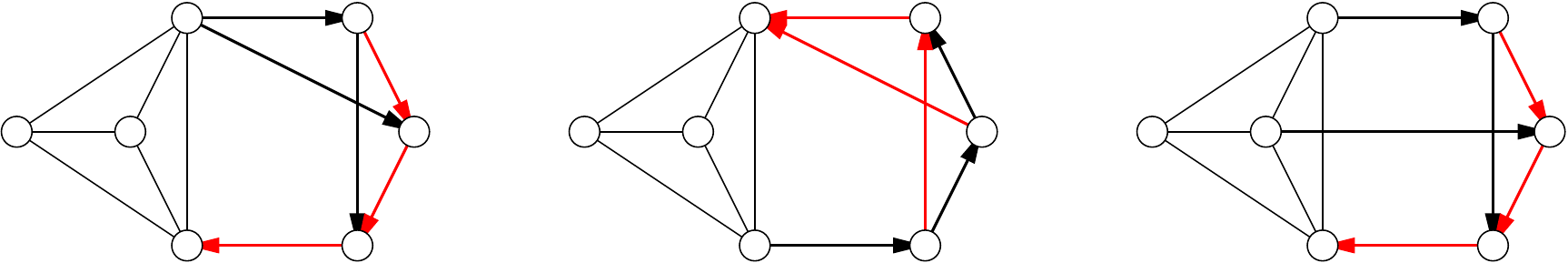}}
\caption{How to extend any $(s,t)$-triple $(\leq{},I,O)$ of $K_4$ to the six remaining edges, depending on whether the vertex of degree $5$ is smaller (left/right) or larger (middle/right) than the vertex of degree $4$ according to $\leq$.}\label{fig:f4n}
\end{center}
\end{figure}

\begin{figure}[H]
\begin{center}
\scalebox{0.6}{\includegraphics{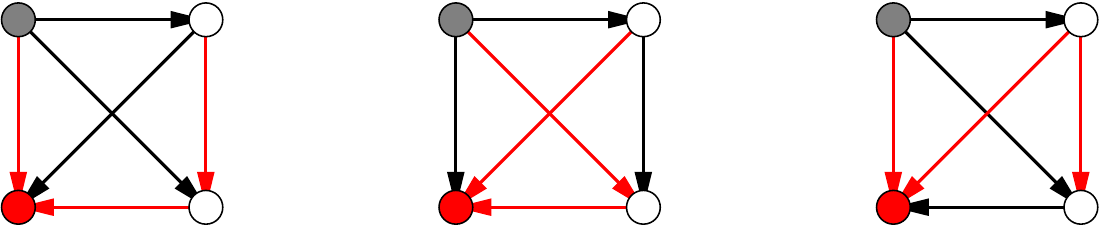}}
\caption{Three $(s,t)$-triples of $K_4$, $s$ and $O$ in black, $t$ and $I$ in red. The two on the left show that we can avoid a specific pair of arcs at some non-root $x$
such that one is the unique in-edge from $O$ at $x$ and the other is the unique out-edge from $I$ at $x$. All three $(s,t)$ triples show that we can
freely choose the unique in-edge from $O$ at $t$ (This  also follows from Theorem \ref{thm:genericgood}).}\label{fig:f5n}
\end{center}
\end{figure}

Consider $s \not= t$ from $V(G)$.
If $s,t$ are both in $V(H)$ we take an $(s,t)$-triple for $H$, and extend it to an $(s,t)$-triple for $G$ according to Figure \ref{fig:f4n},
depending on the order of the two or three neighbors of $X$ in $V(H)$. If $s,t$ are both in $X$ and $u$ is the vertex in $X \setminus \{s,t\}$
then we take one of the two leftmost $(s,t)$-triples depicted in Figure \ref{fig:f5n} for the quotient $Q=K_4$ of $G$ with respect to the partition $\{\{s\},\{t\},\{u\},V(H)\}$
as to achieve that the unique edge in $O$ entering $H$ and the unique edge from $I$ exiting $H$ have distinct endpoints $a,b \in H$ in $G$ (they will automatically be distinct if $G$ is the right graph in Figure \ref{fig:f3n}).
We then insert, as usual, an $(a,b)$-triple of $X$ as to obtain an $(s,t)$-triple for $G$.
Finally, if, without loss of generality, $s \in X$ and $t \in V(H)$, then we take one of the three $(s,V(H))$-triples for $Q$ as to achieve
that the unique edge in $O$ entering $H$ has an endpoint $a \in H$ distinct from $t$. We then extend by an $(a,t)$-triple for $V(H)$
as to obtain an $s,t$-triple for $G$.
\hspace*{\fill}$\Box$

By increasing the degree bound by $1$, we can actually achieve the preconditions of Lemma \ref{lem:L1n}.

\begin{lemma}
  Suppose that $n\in \{4,5,6,7\}$ and $G$ is a simple graph on $n$ vertices with $\delta(G) \geq \lfloor n/2 \rfloor+1$.
  Then $G$ has property (*).
\end{lemma}

{\bf Proof.} For $n=4$ it follows that $G \cong K_4$.
For $n=5$ the graph $G$ is the complement of one of the three graphs on five vertices of maximum degree $1$ and, thus,  has a spanning $4$-wheel,
which certifies (as a generic circuit) that $G$ has property (*).
For $n \in \{6,7\}$ we get $\delta(G) \geq 4$, again look at a set $X \subseteq V(G)$ with $1 \leq |X| \leq 3$, and observe that there
are at least four edges connecting a vertex from $X$ to a vertex from $V(G) \setminus X$; so $G$ is $4$-edge-connected and, thus,
has two edge-disjoint spanning trees. $G$ cannot be the exceptional graph of Lemma \ref{lem:L1n}, as the latter has minimum degree $3$,
so  Lemma \ref{lem:L1n} implies that $G$ has property (*).
\hspace*{\fill}$\Box$

Consider again a simple graph $G$ on $n \geq 4$ vertices with $\delta(G) \geq \lfloor n/2 \rfloor$ and suppose it has a cutvertex $x$.
For each component $C$ of $G-x$, $|V(C)| \geq \lfloor n/2 \rfloor$ holds with equality if and only if $G[V(C) \cup \{x\}] \cong K_{\lfloor n/2 \rfloor +1}$,
As $G-x$ has at least two components $C,D$ we see that $n \geq |V(C)|+|V(D)|+|\{x\}| \geq 2 \lfloor n/2 \rfloor +1$,
which implies that $n$ is odd and $G$ is obtained from identifying two copies of $K_{\lceil n/2 \rceil}$ at a single vertex.
Any simple graph on $n$ vertices with $\delta(G) \geq \lfloor n/2 \rfloor$ nonisomorphic to this exception is $2$-connected.

\begin{theorem}
  \label{thm:dense}
  Suppose that $n \geq 8$ and $G$ is a simple graph on $n$ vertices with $\delta(G) \geq \lfloor n/2 \rfloor$.
  Then $G$ has property (*) unless $n$ is odd and $G$ is obtained from identifying two copies of $K_{\lceil n/2 \rceil}$ at a single vertex.
\end{theorem}

{\bf Proof.} Induction on $n$. We start as in the proof of Lemma \ref{lem:L1n}. We may assume that the given graph
$G$ on $n$ vertices with $\delta(G) \geq \lfloor n/2 \rfloor$ is not the exceptional graph and, hence it is  $2$-connected.
Let $d:=\lfloor n/2 \rfloor \geq 4$. Since $G$ is $d$-edge-connected, it contains a spanning 2T-graph and,
hence, a (not necessarily spanning) generic circuit $C$. As $C$ has property (*) by Theorem
\ref{thm:genericgood}, $C$ extends to a maximal subgraph $H$ of $G$ having property (*),
and $|V(H)| \geq |V(C)| \geq 4$. If $H$ spans $G$ then we are done, otherwise $X:=V(G) \setminus V(H)$ is not empty.
Set $G':=G[X]$. Then $n':=|V(G') \leq n-4$.
By maximality of $H$, (**) every vertex from $X$ has at most one neighbor in $V(H)$, 
so that $\delta(G') \geq d-1= \lfloor n/2 \rfloor -1 \geq \lfloor (n'+4)/2 \rfloor -1 = \lfloor n'/2  \rfloor +1$.
From the first inequality in that chain we get $n' \geq d \geq 4$.
If $n' \leq 7$ then $G'$ has property (*) by Lemma \ref{lem:L1n}. Otherwise, induction applies and shows that $G'$ has property (*) (observe
that the exceptional graph for $n'$ has minimum degree $\lfloor n'/2 \rfloor$, so $G'$ is not exceptional).
Since $G$ is $2$-connected, there are two disjoint edges $ac,bd$ with $a,b \in V(G')$ and $c,d \in V(H)$,
implying that the spanning subgraph of $G$ obtained from the union of $G'$ and $H$ by adding the edges $ac,bd$
has property (*), too, contradicting the maximality of $H$.
\hspace*{\fill}$\Box$

As every complete graph on at least $k\geq 4$ vertices contains a spanning wheel $W_k$ and wheels are generic circuits, we get the following.

\begin{corollary}
  Every graph $G$ on $n\geq 8$ vertices and minimum degree $\delta{}(G)\geq \lfloor{}n/2\rfloor$ has a good $(s,t)$-triple for least $(\lfloor{}n/2\rfloor)^2$ distinct pair of vertices $s,t$. 
  
  \end{corollary}

\section{Remarks and open problems}

\begin{problem}
  \label{prb:goodtransits}
  Does every quartic have an $(s,t)$-triple for some pair $s,t$ of transits?
\end{problem}

\begin{figure}[H]
\begin{center}
\scalebox{0.5}{\includegraphics{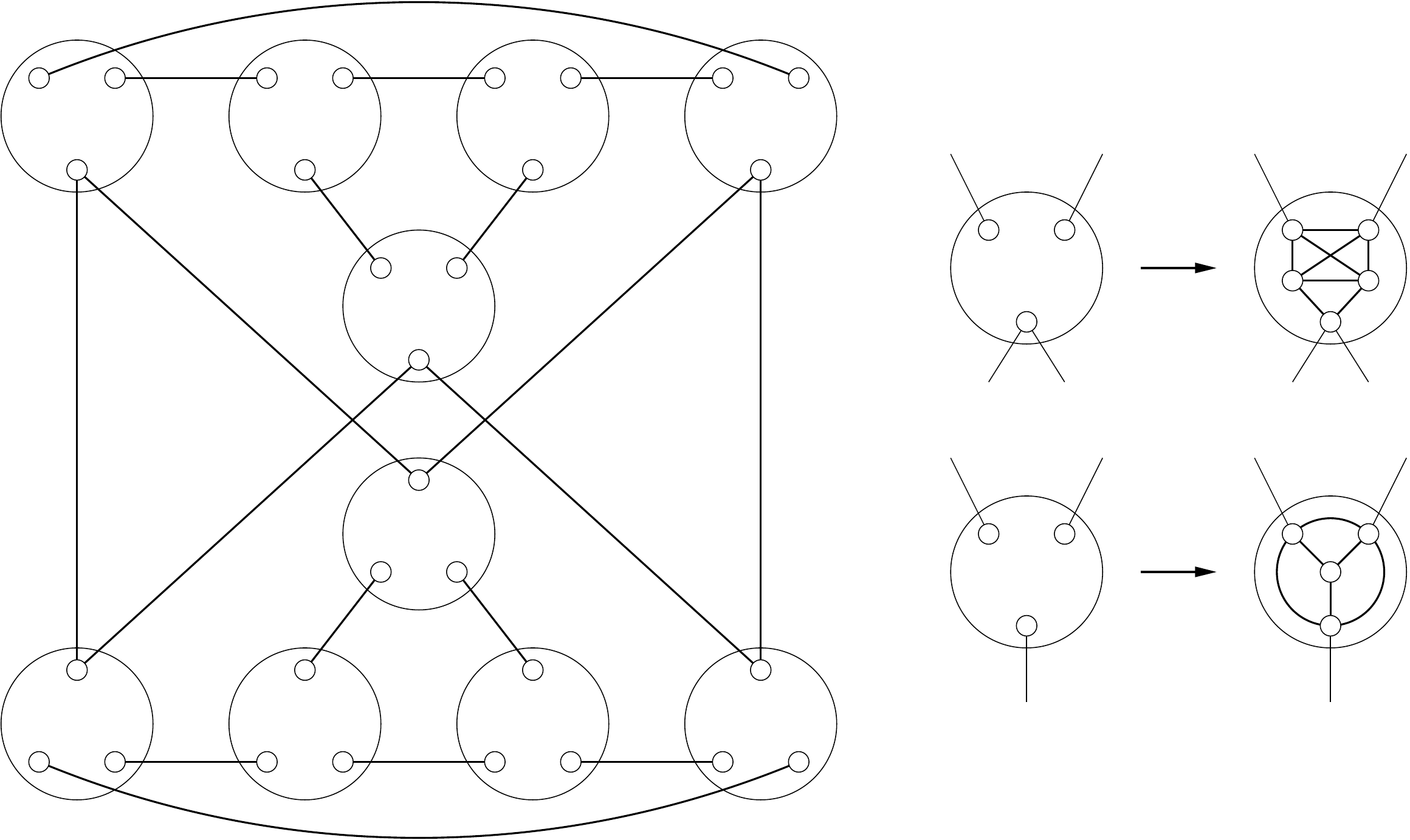}}
\caption{The hypergraph $H$ and how to replace its $3$-edges.}
\end{center}
\end{figure}\label{fig:f5}

Although we do not know a single quartic which does not admit an $(s,t)$-triple for {\em some} of its transits,
we can construct one such that there are pairs of transits for which there is no  $(s,t)$-triple. It not normal showing that
the normality assumption in Theorem \ref{thm:normalisgood} cannot be omitted.

In fact, its four transits can be partitioned into two sets of size $2$ such that
there is an $(s,t)$-triple only if $s,t$ come from distinct sets. We start with the hypergraph $H$ on the left hand side in Figure \ref{fig:f5}.
There are two types of $3$-hyperedges, and we substitute each $e \in E(H)$ with a graph $S_e$ with $V(e) \subseteq V(S_e)$   
according to the rules depicted on the right hand side of Figure \ref{fig:f5}.
Since all substitutes as well as the quotient graph $Q$ of $H$ with respect to its partition into $3$-edges are 2T-graphs (in fact, $Q$ is
a generic circuit), we get a 2T-graph $G$ this way. One readily checks that $G$ is a $3$-connected quartic with $46$ vertices
whose four transits belong to the four substitutes for those three $3$-hyperedges incident with only three edges.

Now suppose that there is an $(s,t)$-triple $T=(\leq,I,O)$ in $G$.
For any induced 2T-subgraph $S$ of $G$, and in particular, for the ten substitutes,
the restriction $(\leq_{|V(S) \times V(S)},I[V(S)],O[V(S)])$ of $T$ to $S$ is an $(a,b)$-triple, where $a$ and $b$ are the minimum and
maximum elements of $S$, respectively. We call these the {\bf local out-root} and {\bf local in-root} of $S$, respectively.
If $a \not= s$ then there exists an edge in $O$ from some vertex $v$ in $V(G) \setminus S$ to $a \in S$,
and this is the unique edge in $O$ from $V(G) \setminus S$ to $S$. Likewise, if $b \not= t$, then there exists a
unique edge in $I$  from $b$ to some vertex $w$ in $V(G) \setminus S$. We say that $v$ {\bf supports} $a$, and $w$ {\bf supports} $b$.

Now we construct an auxilary orientation $D$ of $Q$ as follows:
Whenever $e,f \in V(Q)$ are such that a local root $a$ of $S_e$ is supported by some $v$ in $V(S_f)$, we orient
the edge formed by $a,v$ in $Q$ from $a$ to $v$.
If $v$ happens to be a local root of $S_f$ then, clearly, $a$ cannot support $v$ (in particular, no edge is forced to be ``oriented in both directions''),
so that either $v$ is one of the global roots $s,t$ or is supported by some vertex $w \not= v$.
It follows  that if there is an edge in $D$ from $a$ to $v \not\in \{s,t\}$ such that
$v$ is incident with only one $2$-edge in $H$ then $v$ cannot be a local root of $S_f$.
Hence if $S_f$ does not contain $s$ or $t$, then the two vertices distinct from $e \setminus \{v\}$ are local roots of $S_f$.

Let us denote the four $3$-hyperedges in the first row in the drawing of $H$ by $e,f,g,h$ (from left to right), and
let $j$ denote the $3$-hyperedge immediately below them. Let $U$ be the union of the vertex sets of the five corresponding substitutes
$S_e,S_f,S_g,S_h,S_j$ and assume, to the contrary, that $s,t \not \in U$. This implies that each vertex of $U$ will have out-degree 2 in $D$.

\begin{figure}[H]
\begin{center}
\scalebox{0.5}{\includegraphics{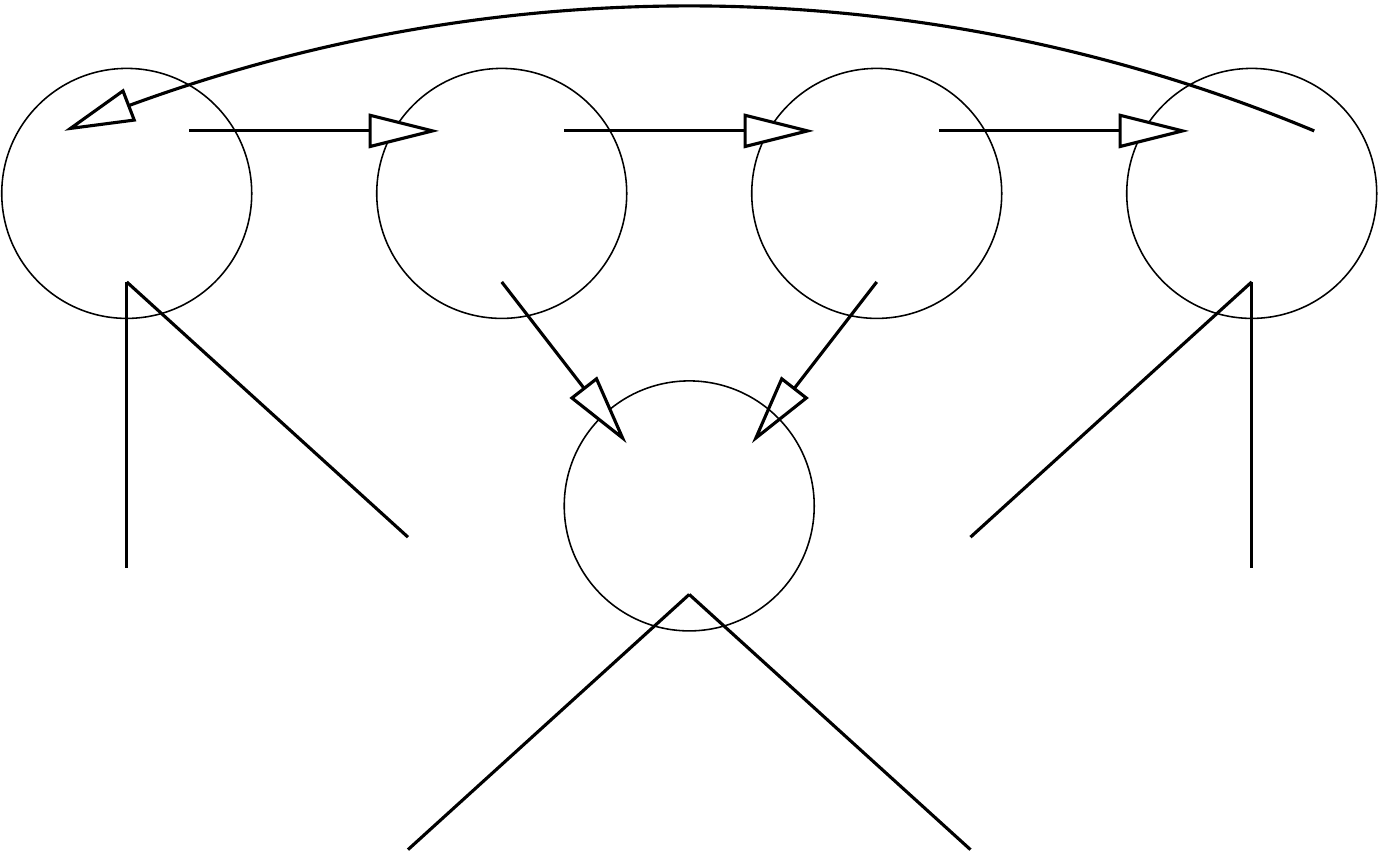}}
\caption{Part of the auxilary digraph $D$. The $4$-cycle $fghe$ forces the arcs $fj,gj$, so that $j$ can have only one local root.}\label{fig:f1n}
\end{center}
\end{figure}

Thus at least one of $fe$ and $fg$ is an arc of $D$, and by the arguments of the previous paragraph we now infer
that either $fehg$ or $fghe$ is a (continuously directed) $4$-cycle in $D$ and  that $D$ contains the arcs $fj,gj$. But then 
 $j$ has only one local root, contradiction. Thus $U$ contains one of $s,t$, and, by symmetry, so does $V(G) \setminus U$.

\begin{figure}[H]
\begin{center}
\scalebox{0.6}{\includegraphics{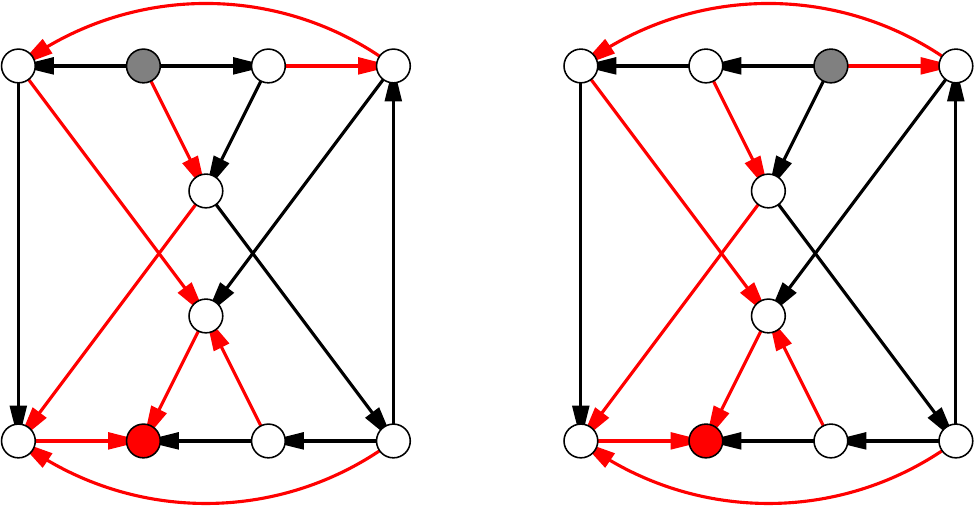}}
\caption{Two $(s,t)$-triples of the quotient quotient $R$ where $s \in P$ (in black) and $t \in Q$ (in red). They can be extended to the edges of the substitutes, as
these allow good triples between all pairs of vertices.}\label{fig:f2n}
\end{center}
\end{figure}

Let $P,Q$ denote the set of the two transits of $G$ in $U$, $V(G) \setminus U$, respectively.
It follows that $G$ is a $3$-connected quartic with the property that whenever there exists an
$(s,t)$-triple
for two transits $s,t$ then one of $s,t$ them is in $P$ and the other one in $Q$. Figure \ref{fig:f2n}
shows that all the 4 pairs of $(s,t)$-triples where $s\in P,t\in Q$ are realizable.\\

If the answer to this Problem \ref{prb:goodtransits} is negative, then we can even prove the following.  
\begin{theorem}
  \label{thm:nogoodor}
  Suppose that  there exists a quartic $Q$ with no  $(s,t)$-triple where $s,t$ are distinct transits. Then the following holds:
  \begin{enumerate}
  \item[(i)] There exists a quartic $G$ with no $(s,t)$-triple for any choice of $s,t\in V(G)$.
  \item[(ii)] There exists a 4-regular  graph $G$ which has no  $(s,t)$-triple for any choice of $s,t\in V(G)$. If $Q$ is 3-connected, then so is $G$ and this would imply that the connectivity bound in Theorem
    \ref{thm:4reg4con} is best possible.
    \end{enumerate}
  \end{theorem}

  \pf Suppose $Q$ is a quartic with no  $(s,t)$-triple where $s,t$ are distinct transits and let $a,b,c,d$ be its transits.\\

  To prove (i) let $G$ be the graph that one obtains from three disjoint copies $Q_1,Q_2,Q_3$ of $Q$ and
  four extra vertices $v_a,v_b,v_c,v_d$ by adding  arcs from $v_x$ to the transit named $x$ in each of $Q_1,Q_2,Q_3$, where $x\in \{a,b,c,d\}$. Clearly $G$ is a quartic and its quotient graph is the generic circuit $K_{3,4}$. Suppose that $G$ has an $(s,t)$-triple $(\leq{},I,O)$ for some choice of its vertices. Then at least one of the copies of $Q$, say w.l.o.g, $Q_1$, does not contain any of the vertices $s,t$. As $Q_1$ is a 2T-graph the restriction of  $(\leq{},I,O)$ to $V(Q_1)$ is an $(s',t')$-triple $(\leq{}',I',O')$ of $Q_1$ where $s',t'$ are transits, contradicting that  $Q_1$ is assumed to have no such pair.\\

  To prove (ii) we construct the graph $W$ from five copies $Q_0,\ldots{},Q_4$ of $Q$ (where we name the transits of $Q_i$ by $a_i,b_i,c_i,d_i$ for $i=0,1,2,3,4$) by adding the edges  $a_ic_{i+1},b_id_{i+1}$ for $i=0,1,2,3,4$, where the indices are taken modulo 5. It is easy to check that $W$ is 4-regular and furthermore, if $Q$ is  3-connected, then so is $W$. Every $(x,y)$-triple of $W$ is  a spanning 2T-graph of $W$ so it is obtained by deleting the edges of a matching of size two from $W$.
  This will leave at least one copy of $Q$ untouched, so as above, assuming that $W$ has an
  $(x,y)$-triple for some pair of vertices $x,y$ we obtain a contradiction to the assumption that $Q$ has no  $(s,t)$-triple where $s,t$ are transits. \qed

  \vspace{2mm}

  In \cite{bangGOpaper} it was shown that a graph $G$ may have arbitrary high edge-connectivity and still have no good orientation when $G$ may have cut vertices.

  \begin{problem}
    Does there exist a natural number $L$ so that every $2$-connected graph which is also $L$-edge-connected has an $(s,t)$-triple for some pair $s,t$ of distinct vertices?
\end{problem}

\begin{problem}
  Does there exists a natural number $K$ so that every $k$-connected graph have an $(s,t)$-triple for some pair $s,t$ of distinct vertices?
\end{problem}

Inspired by Theorems \ref{thm:normalisgood} and \ref{thm:orientalg} the following problem may be doable.

\begin{problem}
  Is there a polynomial algorithm for deciding, for a given quartic $G$ and given distinct
  transits $s,t$ of $G$, whether $G$ has an $(s,t)$-triple?
\end{problem}

\bibliography{refs}
\end{document}

%% file: yaf7.pdf_t
\begin{picture}(0,0)%
\includegraphics{yaf7.pdf}%
\end{picture}%
\setlength{\unitlength}{3947sp}%
\begingroup\makeatletter\ifx\SetFigFont\undefined%
\gdef\SetFigFont#1#2#3#4#5{%
  \reset@font\fontsize{#1}{#2pt}%
  \fontfamily{#3}\fontseries{#4}\fontshape{#5}%
  \selectfont}%
\fi\endgroup%
\begin{picture}(8578,4596)(812,-4350)
\put(7651,-3761){\makebox(0,0)[b]{\smash{{\SetFigFont{12}{14.4}{\rmdefault}{\mddefault}{\updefault}{\color[rgb]{0,0,0}$K_4$}%
}}}}
\put(5701,-1361){\makebox(0,0)[b]{\smash{{\SetFigFont{12}{14.4}{\rmdefault}{\mddefault}{\updefault}{\color[rgb]{0,0,0}$+$}%
}}}}
\put(6901,-2561){\makebox(0,0)[b]{\smash{{\SetFigFont{12}{14.4}{\rmdefault}{\mddefault}{\updefault}{\color[rgb]{0,0,0}$+$}%
}}}}
\put(1351,-3761){\makebox(0,0)[b]{\smash{{\SetFigFont{12}{14.4}{\rmdefault}{\mddefault}{\updefault}{\color[rgb]{0,0,0}$K_4$}%
}}}}
\put(2851,-3761){\makebox(0,0)[b]{\smash{{\SetFigFont{12}{14.4}{\rmdefault}{\mddefault}{\updefault}{\color[rgb]{0,0,0}$K_4$}%
}}}}
\put(4501,-3761){\makebox(0,0)[b]{\smash{{\SetFigFont{12}{14.4}{\rmdefault}{\mddefault}{\updefault}{\color[rgb]{0,0,0}$W_4$}%
}}}}
\put(6151,-3761){\makebox(0,0)[b]{\smash{{\SetFigFont{12}{14.4}{\rmdefault}{\mddefault}{\updefault}{\color[rgb]{0,0,0}$K_4$}%
}}}}
\put(5101,-161){\makebox(0,0)[b]{\smash{{\SetFigFont{12}{14.4}{\rmdefault}{\mddefault}{\updefault}{\color[rgb]{0,0,0}$C$}%
}}}}
\end{picture}%

%% file: go4rdec10.bbl
\begin{thebibliography}{10}

\bibitem{bangJCT51}
J.~Bang-Jensen.
\newblock {Edge-disjoint in- and out-branchings in tournaments and related path
  problems}.
\newblock {\em J. Combin. Theory Ser. B}, 51(1):1--23, 1991.

\bibitem{bangGOpaper}
J.~Bang-Jensen, S.~Bessy, J.~Huang, and M.~Kriesell.
\newblock Good orientations of unions of edge-disjoint spanning trees.
\newblock {\em Submitted}, 2019.

\bibitem{bang2009}
J.~Bang-Jensen and G.~Gutin.
\newblock {\em {Digraphs: Theory, Algorithms and Applications}}.
\newblock Springer-Verlag, London, 2nd edition, 2009.

\bibitem{bangJGT42}
J.~Bang-Jensen, S.~Thomass{\'e}, and A.~Yeo.
\newblock {Small degree out-branchings}.
\newblock {\em J. Graph Theory}, 42(4):297--307, 2003.

\bibitem{bergJCT88}
A.~R. Berg and T.~Jord{\'{a}}n.
\newblock {A proof of Connelly's conjecture on 3-connected circuits of the
  rigidity matroid}.
\newblock {\em J. Comb. Theory, Ser. {B}}, 88(1):77--97, 2003.

\bibitem{edmonds1973}
J.~Edmonds.
\newblock {Edge-disjoint branchings}.
\newblock In {\em {Combinatorial Algorithms}}, pages 91--96. Academic Press,
  1973.

\bibitem{lamanJEM4}
G.~Laman.
\newblock On graphs and rigidity of plane skeletal structures.
\newblock {\em J. Eng. Math.}, 4:331–--340, 1970.

\bibitem{lovaszJCT21}
L.~Lov{\'a}sz.
\newblock {On two min--max theorems in graph theory}.
\newblock {\em J. Combin. Theory Ser. B}, 21:96--103, 1976.

\bibitem{nashwilliamsJLMS39}
C.St.J.A. Nash-Williams.
\newblock {Decomposition of finite graphs into forests}.
\newblock {\em J. London Math. Soc.}, 39:12, 1964.

\bibitem{recski1989}
A.~Recski.
\newblock {\em {Matroid theory and its applications in electric network theory
  and in statics}}.
\newblock Springer-Verlag, Berlin, 1989.

\bibitem{tutteJLMS36}
W.T. Tutte.
\newblock {On the problem of decomposing a graph into {$n$} connected factors}.
\newblock {\em J. London Math. Soc.}, 36:221--230, 1961.

\end{thebibliography}
